\documentclass[10pt]{article}
\usepackage{times}

\usepackage{latexsym}
\usepackage{amsmath}
\usepackage{amssymb}

\newtheorem{Theor}{Theorem}[section]

\newtheorem{Lemma}[Theor]{Lemma}

\def\Dem{\hbox{\textsc {Proof.}\,\,}}
 
\def\Id{\hbox{\sffamily {Id}\,}}

\def\one{\hbox{\sffamily {1}\,}}

\def\Extr{\hbox{\sffamily {Extr}\,}}

\def\One{\hbox{\bf 1}\,}

\def\qed{\hfill $\Box $}

\begin{document}

\title{Hilbert C$^*$-modules are JB$^*$-triples.}

\author{Jos\'e M. Isidro \thanks{Supported by Ministerio de 
Educaci\'on y Cultura of Spain, Research Project PB 98-1371.}\\
Facultad de Matem\'aticas, \\
Universidad de Santiago,\\
Santiago de Compostela, Spain.\\
{\normalsize\tt jmisidro@zmat.usc.es}
}

\date{December 1,  2001}

\maketitle

\begin{abstract}
  We show that every Hilbert C$^*$-module $E$ is a 
		JB$^*$-triple in a canonical way and establish 
		an explicit expression for the holomorphic 
		automorphisms of the unit
		ball of $E$. 
\end{abstract}


\subsection{Introduction} 

Hilbert C$^*$-modules first appeared in 1953 in a work of 
Kaplansky \cite{KAPL} who worked only with modules 
over commutative unital C$^*$-algebras. In 1973 
Paschke \cite{PAS} proved that most of the properties of Hilbert 
C$^*$-modules were valid for modules over an arbitrary C$^*$-algebra. 
About the same time Reiffel independently developed much of the 
same theory and used it to study representations of C$^*$-algebras. 
Since then the subject has grown and spread rapidly and now there is an 
extensive literature on the topic (see \cite{LAN2} for a 
systematic introduction). Many interesting developments have been made 
by Kasparov, who used Hilbert C$^*$-modules as the framework 
for K-theory. More recently Hilbert C$^*$-modules have been a 
useful tool in the C$^*$-algebraic approach to quantum groups. 
The geometry of Hilbert C$^*$-modules has been investigated 
by Solel in \cite{SOL}, where the isometries 
of these Banach spaces have been characterized. See also \cite{AND}.

On the other hand Kaup, searching for a
metric-algebraic setting in which he could make the study of bounded
symmetric domains in complex Banach spaces, introduced a class of complex
Banach spaces called  JB$^*$-triples. In 1983 he proved 
that, except for a biholomorphic bijection, every such a domain is the 
open unit ball of a JB$^*$-triple \cite{KAUR}. In 1981 he
made the complete analytic classification of bonded symmetric domains in
reflexive Banach spaces \cite{KAUS}. Since then the study of JB$^*$-triples
has grown and spread considerably.  

These two theories have developed independently from one another. Here 
we show that every Hilbert C$^*$-module is, in a canonical way, a
JB$^*$-triple, a bridge between the two theories that may be useful in
the study of the geometry of Hilbert C$^*$-modules. We also establish an
explicit expression of the  holomorphic automorphisms of the unit ball of a
Hilbert C$^*$-module.


\subsection{Hilbert C$^*$-modules}

We now introduce formally the objects we shall be studying. 
Let $A$ be a C$^*$-algebra (not necessarily unital or commutative) 
where the product is denoted by juxtaposition $xy$, the norm is $\Vert \cdot
\Vert_A$ and the $x\mapsto x^*$ stands for the conjugation. An {\sl inner
product} $A$-module is a complex linear space 
$E$ with two laws of composition $E\times A\to E$ (denoted by 
$(x,a)\mapsto x.a$) and $E\times E\to A$ (denoted by $(x,y)\mapsto 
\langle x ,y\rangle$) such that the following properties hold: 
\begin{enumerate}
\item With respect to the operation $(x,a)\mapsto x.a$, $E$ is a right
$A$-module with a compatible scalar multiplication, that is, $\lambda (x.a) 
=(\lambda x).a= x.(\lambda a)$ for all $x\in E$, $a\in A$ and $\lambda \in
\mathbb{C}$. 
\item The inner product $(x,y)\mapsto\langle x,y\rangle$ satisfies 
\begin{eqnarray}
\langle x, \alpha y+\beta z\rangle&=&\alpha\langle x,y\rangle + \beta \langle
x,z\rangle, \\
\langle x, y.a\rangle&=&\langle x,y\rangle a\\
\langle y,x\rangle&=&\langle x,y\rangle^* \\
\langle x,x\rangle &\leq &0, \qquad \hbox{\rm if }\;\; \langle x, x\rangle =0 
\;\; \hbox{\rm then}\;\; x=0.
\end{eqnarray}
for all $x,yz\in E,$ all $\alpha ,\beta \in \mathbb{C}$ and all $a\in A$.  
\end{enumerate}
Note that in particular, the inner product is complex linear in the second 
variable while it is conjugate linear in the first. This convention is in line
with the recent research literature. Let $E$ be an inner product $A$-module; 
then the Cauchy-Schwarz inequality 
$$
\langle y,x\rangle\,\langle x,y\rangle \leq \Vert \langle x,x\rangle\Vert \,
\langle y,y\rangle,
$$
holds, hence $\Vert x\Vert_E ^2\colon = \Vert \langle x,x\rangle\Vert_A $
is a norm in $E$ with respect to which the inner product and
the module product are continuous, that is 
$$ 
\Vert \langle x,y\rangle\Vert _A\leq \Vert x\Vert_E\,\Vert y\Vert_E, 
\qquad \Vert x.a\Vert \leq \Vert x\Vert_E\, \Vert a\Vert_A. 
$$ 
To simplify the notation we shall use the same symbol $\Vert \cdot \Vert$ to
denote the norms on $A$ and $E$. We can also define an $A$-valued {\sl norm} 
by $\vert x\vert\colon = \langle x,x\rangle^{\frac{1}{2}}$ for $x\in E$ and 
we have 
$$
\vert \langle x,y\rangle\vert \leq \Vert x\Vert \, \vert y\vert , \qquad 
\vert \langle x,y\rangle \vert \leq \vert x\vert \,\Vert y\Vert, 
$$
and the module product is also continuous with respect to the new norm on $A$ 
since 
$$
\Vert x.a\Vert \leq \Vert x\Vert \, \vert a\vert.
$$
However the $A$-valued norm in an inner $A$-module $E$ need to be handled with
care. For example it need not be the case $\vert x+y\vert \leq \vert x\vert
+\vert y\vert $. 
An inner product $A$-module $E$ which is 
a Banach space with respect to the norm $\Vert \cdot\Vert $ is called a
{\sl Hilbert C$^*$-algebra  module}. Every C$^*$-algebra can be converted into 
a Hilbert C$^*$-algebra module by taking $E\colon = A$ with the natural 
module operation $x.a\colon =xa$ and the inner product $\langle a,b\rangle
\colon = a^*\,b$ for $x,a,b\in A$. 

A linear map $f\colon E\to E$ is called an $A$-map if $f(x.a)=f(x).a$ holds for 
all $x\in E$ and $a\in A$, and we say that $f$ is {\sl adjointable} if there 
exists an $A$-map $f^*\colon E\to E$ such that 
$$
\langle f(x),y\rangle =\langle x, f^*(y)\rangle , \qquad x,y\in E.
$$
In such a case $f$ is continuous (though the converse is not true!), $f^*$ is
adjointable and
$(f^*)^*=f$.  We let $\mathcal{A}(E)\subset \mathcal{L}(E)$ denote the vector  
space of all adjointable $A$-module maps on $E$. In fact $\mathcal{A}(E)$ is a 
C$^*$-algebra in the operator norm since $\Vert f^*f\Vert = \Vert f\Vert^2$ 
holds for all $f\in \mathcal{A}(E)$. For $x,\,y\in E$ we define $\theta_{x,y}$
(also denoted by $x\otimes y^*$) by 
$$ \theta _{x,y}(z)\colon = x.\,\langle y,z\rangle, \qquad z\in E,$$ 
Then $\theta_{x,y}$ is adjointable and $\theta_{x,y}^*=\theta_{y,x}$ 
(see \cite{LAN2} p. 9).  
For later reference we state the following 
\begin{Lemma}\label{pos}
Let $E$ be a Hilbert C$^*$-module and let $f\colon E\to E$ be a bounded 
$A$-module map. Then  
$f$ is a positive element in the C$^*$-algebra $\mathcal{A}(E)$ if and only if 
$\langle x,f(x)\rangle\geq o$ for all $x\in E$.
\end{Lemma}
We refer to \cite{LAN2} for background on Hilbert C$^*$-modules and for the 
proofs of the above results.


\subsection{JB$^*$-triples.}

For a complex Banach space $X$ denote by $\mathcal{L}(X)$ 
the Banach algebra of
all bounded complex-linear operators on
$X$. A complex Banach space $Z$ with a 
continuous mapping $(a, b, c)  
\mapsto \{a, b, c\}$ from $Z\times Z\times Z$ to $Z$ is called a {\it JB*-triple} 
if the following conditions are satisfied for all $a, b, c, d \in Z$, where 
the operator $a\square b\in \mathcal{L}(Z)$ is defined by $z\mapsto \{abz\}$ and 
$\lbrack\, , \, \rbrack$ is the commutator product:
\begin{enumerate}
\item $\{abc\}$ is symmetric complex linear in $a, c$ and conjugate linear 
in $b$.
\item $\lbrack a\square b , \, c\square d \rbrack = \{a,b,c\}\square d -
c\square \{d,a,b\}$ (called the Jordan identity. ) 
\item $a\square a$ is hermitian and has spectrum $\geq 0.$
\item $\Vert \{a,a,a\}\Vert = \Vert a\Vert ^3$. 
\end{enumerate}
If a complex vector space $Z$ admits a JB*-triple structure, then the norm and 
the triple product determine each other. 
An {\it automorphism} is a linear bijection $\phi \in \mathcal{L}(Z)$ such 
that $\phi \{z,z, z\}= \{(\phi z), ( \phi z),(\phi z)\}$ for $z\in Z$, which 
occurs if and only if $\phi$ is a surjective linear isometry of $Z$.

Recall that every C*-algebra $Z$ is a JB*-triple with respect
to the triple product $2\{abc\} \colon =(ab^*c+cb^*a)$. In that case, 
every  projection in $Z$ is a tripotent and more generally the tripotents
are precisely the partial isometries in $Z$. C$^*$-algebra derivations
and C$^*$-automorphisms are derivations and automorphisms of $Z$ as a 
JB$^*$-triple though the converse is not true. 

We refer to \cite{KAUR}, 
\cite{UPM} and the references therein for
the background of JB$^*$-triples theory.


\subsection{Hilbert $C^*$-modules are JB$^*$-triples.}
For $a\in A$ fixed, we denote by $R_a \in
\mathcal{L}(E)$ the operator $x\mapsto x.a$ of right
multiplication by $a$. 
\begin{Theor}
Every Hilbert C$^*$-module $E$ is a a
JB$^*$-triple in a canonical way
\end{Theor}
\Dem Let $E$ be a Hilbert C$^*$-module over
the C$^*$-algebra $A$ and define a triple 
product in $E$ by 
\begin{equation}\label{tp}
2\,\{x,y,z\}\colon = x.\langle
y,z\rangle+z.\langle y,x\rangle, \qquad
x,y,z\in E.
\end{equation} 
It is clear that $\{\cdot , \cdot , \cdot \}$
symmetric complex linear in the external variables, 
and complex conjugate linear in the middle variable. 
It is a matter of routine calculation to check that
the triple product satisfies the Jordan identity. On
the other hand, for fixed $x\in E$ we have 
$$2 \,(x\Box x)\, z= x.\langle x, z\rangle +
z.\langle x,x\rangle , \qquad z\in E
$$ 
which can be written in the form $x\Box x=\frac{1}{2} 
(\theta_{x,x}+R_{\vert x\vert ^2})$. We show that the  
summands in the right hand side of the latter are
hermitian elements in the algebra $\mathcal{L}
(E)$. Since $\mathcal{A}(E)$ is a closed complex
subalgebra of $\mathcal{L}(E)$ and contains the
unit element, it suffices to consider the numerical
range of $\theta_{x,x}$ and $R_{\vert x\vert ^2}$
viewed as elements in the C$^*$-algebra
$\mathcal{A}(E)$, and we have seen before that 
$\theta_{x,x}$ is selfadjoint. Clearly  
$$
(\exp \,itR_{\vert x\vert ^2})\,(w)= w.(\exp \,it \vert
x\vert^2) \qquad w\in E,
$$
and as $\exp \,it \vert x\vert^2$ is a unitary element
in $A$, the operator $\exp \,itR_{\vert x\vert
^2}$ is an isometry of $E$ for all $t\in \mathbb{R}$, 
which shows that $R_{\vert x\vert ^2}$ is hermitian. 
For $y\in E$ we have 
$$
\langle y,\; \theta_{x,x}(y)\rangle= \langle y,\;
x.\langle
x,y\rangle\rangle=\langle y,x\rangle\langle
x,y\rangle\geq 0$$  which by (\ref{pos}) proves that
$\theta_{x,x}\geq 0$ in $\mathcal{A}(E)$ hence also in
$\mathcal{L}(E)$.  Clearly $\vert x\vert^2\geq 0$ in
$A$, hence its spectrum satisfies $\sigma_A(\vert
x\vert^2)\subset [0,\infty )$ and therefore 
$$
\sigma_{\mathcal{L}(A)}(R_{\vert x\vert^2})\subset 
\sigma_A(\vert x\vert^2)\subset [0,\infty )
$$
Since the numerical range is the convex hull of the
spectrum, $R_{\vert x\vert^2}\geq 0$ as we wanted to
check. 

Let us set $y\colon =\langle x,x\rangle \in
A$ for every $x\in E$. The definition of the norm in 
$E$ and the properties of the norm in the
C$^*$-algebra $A$ yield 
\begin{eqnarray*}
\Vert \{ x,x,x\}\Vert ^2&=&\Vert x.\langle
x,x\rangle\Vert ^2=
\Vert \langle x.\langle x,x\rangle ,\; 
  x.\langle x,x\rangle \rangle \Vert =\\
\Vert \langle x,x\rangle \,\langle x,x\rangle\,\langle
x,x\rangle \Vert &=& \Vert \{y,y,y\}\Vert = \Vert
y\Vert ^3 = \Vert \langle x,x\rangle \Vert ^3 = \Vert
x\Vert ^6
\end{eqnarray*}
which shows property 4. Finally, this is the unique
JB$^*$-triple structure on $E$ since the triple
product is determined by the norm of $E$. \qed


\subsection{Holomorphic automorphisms of the
unit ball.}

Motivated by the deep formal analogy between Hilbert
C$^*$-modules $E$ and Hilbert spaces $H$, we shall
establish an explicit formula for the holomorphic
automorphisms of the unit ball of $E$. 
Recall \cite{KAUR} that, for $c\in E$, the
{\sl Bergmann operator} of $E$ is given by 
$$ 
B(c,c)(x)\colon = x-2(c\Box c)(x)+ Q_c^2(x), \qquad
x\in E.
$$ 
In our case 
\begin{eqnarray*}
2(c\Box c)(x) &=& 2\{c,c,x\}=c.\langle
x,x\rangle +x.\vert c\vert^2= c\otimes c^*
(x)+x.\vert c\vert^2, 
\\ 
 Q_c^2(x)&=& \{c ,\,Q_c(x),\,c\}= \{c,\;
c.\langle x,c\rangle ,\, c\}= \\ 
c.\langle c.\langle x,c\rangle ,
\,c\rangle &=& 
c.\langle c,x\rangle \vert c\vert ^2=
(c\otimes c^*) (x.\vert c\vert^2).
\end{eqnarray*}
Therefore 
\begin{eqnarray*}
B(c,c)(x)&=& x.(\One -\vert c\vert ^2 )+ 
(c\otimes c^*)\,\big( x.(\One -\vert c\vert ^2 )\big)
=\\ 
{}&=& ( \One -c\otimes c^* ) \, \big( 
x.(\One -\vert c\vert ^2 ) \big) . 
\end{eqnarray*} 
Recall that $\One -c\otimes c^*$ and $\one -\vert
c\vert ^2$ are selfadjoint
elements in the C$^*$-algebras $\mathcal{A}(E)$ and
$A$ respectively, hence they have well defined square
roots. We show the operator  
$$
B_c(x)\colon = (\One -c\otimes c^*)^{\frac{1}{2}} 
\,\big( x. (\one -\vert c\vert
^2)^{\frac{1}{2}}\big), 
\qquad x\in E,  
$$ 
satisfies $B_c^2 =B(c,c)$. Indeed, since $\One
-c\otimes c^*$ is an $A$-linear map so is
its square root and we have 
\begin{eqnarray*}
&{}&B_c(B_c(x))= (\One -c\otimes c^*)^{\frac{1}{2}}\; 
\big( B_c(x).\, (\one -\vert c\vert^2)^{\frac{1}{2}}   
\big) =  \\
&{}&\big( \,(\One -c\otimes c^*)^{\frac{1}{2}}\,
B_c(x)\,\big) .  
(\one -\vert c\vert^2)^{\frac{1}{2}})= \\
&{}&\big( \,(\One -c\otimes c^*)^{\frac{1}{2}}\,
[\, (\One -c\otimes c^*)^{\frac{1}{2}} \,x. (\one
-\vert c\vert^2)^{\frac{1}{2}}\,]\big)\, .\, 
(\one -\vert c\vert^2)^{\frac{1}{2}} =\\
&{}&(\One -c\otimes c^*) x (\one -\vert c\vert^2)= 
B(c,c)(x) 
\end{eqnarray*}
as we wanted to check.  

Now we prove that for $c$ and $x$ in the open unit
ball of $E$, $\one +\langle c,x\rangle$ is an
invertible element in $A$. Indeed, let us denote by 
$\sigma (a)$ and $v(a)$ the spectrum and the numerical
range of $a$ in the unital algebra $A$. Then 
$$
\sigma (\one +\langle c,x\rangle )\subset v(\one
+\langle c,x\rangle )\subset 1+v(\langle c,x\rangle).
$$
Since $\Vert \langle c.x\rangle \Vert \leq \Vert
c\Vert \, \Vert x\Vert <1$, the numerical range
$v(\langle c,x\rangle )$ is contained in the open
unit disc of $\mathbb{C}$, therefore $-1\notin 
v(\langle c,x\rangle )$ and by the above $0\notin
\sigma (\one +\langle c,x\rangle )$. In particular, 
$$
x.(\one +\langle c,x\rangle )^{-1}
$$
is well defined in $A$. Recall \cite{KAUR} 
that for
$c$ in the open unit ball of $E$, the {\sl
transvection} $g_c$ is the holomorphic automorphism of 
the open ball of $E$ given by 
$$
g_c(x)\colon = c+ B_c \big( x(\One +c\Box
x)^{-1}\big),
\qquad \Vert x\Vert <1. 
$$
Replacing the expression of $B_c$ and $\One+c\Box x$ 
we get
$$
g_c(x) =c+ (\One -c\otimes c^*)^{\frac{1}{2}}\;[\,
x.(\one +\langle c,x\rangle )^{-1}\, (\one -\vert
c\vert ^2)^{\frac{1}{2}}\,]
$$
an expression that can be restated in terms of
projections and coincides with the well
known formula for the transvections of the ball in a
Hilbert space see (\cite{HAR} p. 21). By
\cite{KAUR} every holomorphic automorphism
$h$ of the unit ball of $E$ can be represented
in the form $h= L \circ g_c$ for some
surjective linear isometry of $E$ and some
$c\in E$ with $\Vert c\Vert <1$.

\subsection{Extreme points of the unit ball.}
For a complex Banach spce $E$, the set $\Extr B_E$ of 
extreme points in the unit ball $B_E$ of $E$ palys an 
important role in the study of the geometry of $E$. 
Obviously, we can replace a HIlbert C$^*$-module $E$ 
with its associated JB$^*$-triple in order to study 
the extreme points of the ball $B_E$. By (\cite{KAUP} 
prop. 3.5) we have 
$$
\Extr B_E=\{c\in E \colon B(c,c)=0\}
$$
that is, for $c\in E$ the condition $c\in \Extr B_E$ 
is equivalent to 
$( \One -c\otimes c^* ) \, \big( 
x.(\One -\vert c\vert ^2 ) \big)$ for all $x\in E.$ 
Therefore we have two obvious families of extreme points 
given by 
\begin{eqnarray*}
E.(\one -\vert c\vert) &=&\{0\}\Longrightarrow c\in \Extr B_E,\\
(\Id -c\otimes c^*)\,E&=&\{0\}\Longrightarrow c\in \Extr B_E.
\end{eqnarray*}
These two families may coincide (as it occurs when $E$ is Hilbert
space) but in general they are different. We do not know whter 
every extreme point lies in one of the above families. Every
extreme point is a tripotent, that is, it satisfies $c=\{c, c,c\}= 
c.\langle c, c\rangle$. 

It might be interesting to characterize Hilbert C$^*$-triples 
within the category of JB$^*$-triples.


\end{document}